\documentclass{amsart}
\usepackage{a4}
\usepackage{amssymb,amsmath,amsfonts,amsthm,mathrsfs,dsfont}
\usepackage{graphicx}
\usepackage{color}

\newtheorem{theorem}{Theorem}
\newtheorem{lemma}{Lemma}[section]

\newtheorem{proposition}[theorem]{Proposition}

\theoremstyle{remark}
\newtheorem{remark}[theorem]{Remark}

\def\Carre#1#2{\vbox{
   \hrule height .#2pt
   \hbox{\vrule width .#2pt height #1pt \kern #1pt
      \vrule width .#2pt}
   \hrule height .#2pt}}
\def\GC{\text{G\!C}}
\def\BGC{\text{B\!G\!C}}
\def\e{\varepsilon}
\def\R{\mathbb{R}}

\def\Z{\mathbb{Z}}

\def\H{\mathcal{H}}

\def\LM#1{\hbox{\vrule width.2pt \vbox to#1pt{\vfill \hrule width#1pt
height.2pt}}}
\def\LL{{\mathchoice {\>\LM7\>}{\>\LM7\>}{\,\LM5\,}{\,\LM{3.35}\,}}}
\def\restr{{\LL}}
\def\Div{\textup{div}\,}
\def\car#1{\chi_{#1}}
\def\dist{\textup{dist}}
\def\ds{\displaystyle}

\def\wto{\rightharpoonup}
\def\E{\mathbf{e}}
\def\EE{\mathbf{E}}
\def\Ed{\E_d}
\def\EEd{\EE_d}
\def\ED{E}
\def\Dp{\Div^+}
\def\Dm{\Div^-}
\def\Tr{{\textup{Tr}\,}}
\def\Trp{{\textup{Tr}^+}}
\def\Trm{{\textup{Tr}^-}}
\def\jmp#1{\blue{[#1]}}
\def\ee{\eta_\e}

\newcommand{\blue}[1]{{#1}}

\title[Non interpenetration]{Approximation of a brittle fracture energy with
 a constraint of non-interpenetration}
\author[A. Chambolle] {Antonin Chambolle} 
\address[Antonin Chambolle]{CMAP, Ecole Polytechnique, CNRS,
91128 Palaiseau Cedex, France}
\email[A. Chambolle]{antonin.chambolle@polytechnique.fr}

\author[S. Conti] {Sergio Conti} 
\address[Sergio Conti]{Institut f\"ur Angewandte Mathematik,
Universit\"at Bonn, 53115 Bonn, Germany}
\email[S. Conti]{sergio.conti@uni-bonn.de}

\author[G.A. Francfort]{Gilles A. Francfort} 
\address[Gilles Francfort]{LAGA, Universit\'e Paris-Nord, 
Avenue J.-B. Cl\'ement 93430 - Villetaneuse, France \& Courant Institute, 251 Mercer Street, New York, NY10012, USA}
\email[G. Francfort]{gilles.francfort@univ-paris13.fr}

\begin{document}

\date{\today}

\begin{abstract} Linear fracture mechanics (or at least the initiation part of that theory) can be framed in a variational context as a minimization problem over a $SBD$ type space. The corresponding functional can in turn be approximated in the sense of $\Gamma$-convergence by a sequence of functionals involving a phase field as well as the displacement field.  We show that a similar approximation persists
if additionally imposing a non-interpenetration constraint in the minimization, namely that {only nonnegative normal jumps should be permissible.}

\vskip.5cm



\end{abstract}
\maketitle

%%%%%%%%%%%%%%%%%%%%%%%%%%%%%%%%%%%%%%%%%%%%%%%%%%%%%%%%%%%%%%%%
\section{Introduction}
%%%%%%%%%%%%%%%%%%%%%%%%%%%%%%%%%%%%%%%%%%%%%%%%%%%%%%%%%%%%%%%%
The past twenty years or so have been fertile ground for the development of a variational theory of fracture evolution for brittle materials in the context of globally minimizing energetic evolutions~\cite{FranMari}; see for instance \cite{BFM} for a panorama of the theory as it stood a few years back. One of the key ingredients of that theory is a stability criterion which states that the sum of the elastic energy and of the (add-)surface energy at any given time should be minimal for the actual {(add-)}crack at that time among all (add-)cracks and all compatible displacement fields satisfying the loading requirements at that time;  think for example of a time-dependent boundary displacement on the entire boundary of the domain.

In the case of isotropic linear elasticity\footnote{We will restrict our investigation to  isotropic  and  homogeneous materials, although there are no major obstacles in generalizing the result to arbitrary material symmetries and inhomogeneities, at least as far as the elastic energy is concerned.}, the elastic energy is of the form
$$
\frac{1}{2}\int_\Omega \{2\mu |\EE(u)|^2\,+\,\lambda (\Div u)^2\}\ dx$$
where $u:\Omega\to\R^n,\; n=2,3$ is the displacement field over the domain $\Omega$, $\EE(u):=1/2(\nabla u+\nabla^T u)$ is the linearized strain and
$\mu>0,\lambda>-2\mu/n$ are the Lam\'e {constants} of the elastic material.
%We will assume henceforth that {\it both}  $\mu$ and $\lambda$ are positive which, although not required as far as $\lambda$ is concerned, is  the case for most brittle materials.
In this work it will actually be more convenient
to rewrite this energy, as is classical, in the form
\[
\frac{1}{2}\int_\Omega \{2\mu |\EEd(u)|^2\,+\, K (\Div u)^2\}\ dx
\]
where $K=\lambda+2\mu/n>0$ and $\EEd(u) = \EE(u)-\frac{\Div u}{n} I$ is the
\textit{deviatoric part} of the tensor $\EE(u)$ (that is,
its orthogonal projection onto trace-free tensors).

Following in the footstep of {\sc A.A. Griffith}'s foundational paper  \cite{Gr}, it is further  assumed that the (add-)surface energy is proportional to the surface area of the (add-)\newline crack, the coefficient of proportionality, the toughness,  being denoted by $G_c$.

If  considering an uncracked sample $\Omega$ submitted to a boundary displacement $w$ on its boundary $\partial \Omega$,  the crack initiation 
problem then consists in minimizing the sum of those two contributions among
 all cracks -- say closed sets $\Gamma\subset\bar\Omega$ of finite $\H^{n-1}$-measure -- and all displacements
fields $v$ which lie in $H^1(\Omega'\setminus\Gamma;\R^n)$ with $\Omega\subset\subset \Omega'$ and $u=w$ on $\Omega'\setminus\bar\Omega$, that is
\begin{multline*}
\min\bigg\{ \frac12\int_{\Omega\setminus\Gamma} \{2\mu |\EEd(u)|^2\,+\,K (\Div u)^2\} \,dx+\,G_c\H^{n-1}(\Gamma)\,:
\\ u\in H^1(\Omega'\setminus\Gamma), u\equiv w \text{ on } \Omega'\setminus\bar\Omega \bigg\} \,.
\end{multline*}
We will call this formulation the strong formulation. As was first advocated by  {\sc E. De Giorgi}, it is mathematically
convenient to address the strong formulation  in a weak form as follows
\begin{multline}\label{wf}
\min\bigg\{ \frac12\int_\Omega \{2\mu |\Ed(u)|^2\,+\, K (\Tr \E(u)^2\} \,dx+\,G_c\H^{n-1}(J_u):
\\ u\in GSBD(\Omega'), u\equiv w \text{ on } \Omega'\setminus\bar\Omega.\bigg\}.
\end{multline}
Above, the space $GSBD(\Omega')$ is an adequately defined {variant} of the space $SBD(\Omega)$, the space of special functions with bounded deformations. We refer to e.g. \cite{ACDM,BCDM} for a definition and useful properties of the latter and to \cite{DM} for a definition of the former.
Notationwise, $\E(u)$ denotes the Lebesgue absolutely continuous part of  $\EE(u)$ (which {for  $u\in SBD(\Omega)$} is a bounded Radon measure), $\Ed(u):=\E(u)-\frac{\Tr\E(u)}{n}I$ its deviatoric part, while $J_u$ denotes the jump set of $u$ (see e.g. \cite{ACDM} for a precise definition). In particular, $\Tr \E(u)$, which is the trace of $\E(u)$, is the
absolutely continuous part of the divergence $\Div u$.

Remark that the equivalence between the two formulations is still an open problem. By   contrast, the analogous weak formulation with gradients in lieu of symmetrized gradients  has been shown to be equivalent to its  strong counterpart \cite{CL}. In the linearized elasticity framework, a partial result in this
direction, in dimension~2, was issued last year~\cite{CFIreg},
\blue{and further extended to higher dimension in~\cite{ChambolleContiIurlano}}.

Also remark that the existence of a minimizer for Eq.~\eqref{wf}  remains open, except in the case where $n=2$, thanks to a very recent result  \cite[Theorem 6.1]{FS}.

{From} a computational standpoint a formulation such \blue{as} Eq.~\eqref{wf} is rather useless because the test space for the minimization problem is too singular.  It is also widely acknowledged in various fields of physics that sharp interface models are most profitably, and arguably more realistically, addressed  as limits of phase field type models. In the gradient case, {\sc E. De Giorgi} suggested  an approximating formulation which was later proved to $\Gamma$-converge (in the appropriate topology)  to the sharp interface model in \cite{AT}. That approximation is usually referred to as an {\sc Ambrosio-Tortorelli} type approximation, at least in the mathematical community. 

In the present setting, the  approximating phase field  functional is
\begin{multline*}\label{approxim}
E_\e(u,v) = \frac12\int_\Omega (\ee+v^2)
\left\{2\mu |\EEd(u)|^2+ K (\Div u)^2\right\}
 dx  + G_c\int_\Omega \left\{\e|\nabla v|^2+\frac{(1-v)}{4\e}^2\right\}\ dx
\end{multline*}
with $\ee\ll \e$ and the proof of the $\Gamma$-convergence can be found in \cite{ChaSBD1, ChaSBD2} under the additional constraint that $\|u\|_{L^\infty(\Omega;\R^n)}\le M$ for some fixed constant $M$ (which means in particular that the functional framework can then be restricted to $SBD(\Omega)$).
\blue{This result was extended to the %entire
 space $GSBD(\Omega)$ in \cite{Iurlano}, see also~\cite{CFI17b,ChambolleCrismale2017} for recent developments.}

Our contribution starts with the observation that  the weak formulation of Eq.~\eqref{wf} is unphysical because it fails to account for non-interpenetration, that is for the physically obvious requirement that the crack lips should not interpenetrate.
Such  will be the case at any point $x$ where  $\jmp{u}(x)\cdot\nu(x)<0$ where $\jmp{u}(x)$ denotes the jump at $x$ and $\nu(x)$ the normal to the jump set at $x$, well-defined $\H^{n-1}$-a.e.~on $J_u$.
We will thus require that  $\jmp{u}(x)\cdot\nu(x)\ge0, \; \text{ for $\H^{n-1}$-a.e. } x\in J_u$. This can be viewed as a linearized non-interpenetration condition.\footnote{The reader is directed to \cite{GP,DML,DML2} for a treatment of non-interpenetration in the setting of finite deformations.} The goal of this paper is to establish a result of $\Gamma$-convergence for an approximation of the functional of Eq.~\eqref{wf} \`a la {Ambrosio-Tortorelli}
under that further restriction  on the jump. 

We  observe that our proof will address other types of (convex) constraints,
such as an infinitesimal shear condition~\cite{LancioniRC},
or conditions on the eigenvalues of the strain tensor{~\cite{Freddi-Royer-Carfagni-2011a}}, see Remark~\ref{rem:othermodels} below.

Throughout this paper we assume that $\Omega$ is bounded,
with a boundary which is everywhere locally a continuous graph.
{In the case of the non-interpenetration condition,}
the sharp interface functional is given by
\begin{equation}\label{energy}
E(u)  =\begin{cases}
\ds
\frac12 \int_\Omega \{2\mu |\Ed(u)|^2\,+\, K (\Tr\E(u))^2\} \,dx+\,G_c\H^{n-1}(J_u)
& \textrm{if } u\in SBD(\Omega),
\\ &  \hspace{-1cm}\jmp{u}\cdot\nu\ge 0
\; \H^{n-1}-\textrm{a.e.\ in } J_u
\\
+\infty & \textrm{otherwise.}
\end{cases}
\end{equation}
For $u\in SBD(\Omega)$, the measure $\EE(u)$ decomposes as follows~\cite{ACDM}:
\[
\EE(u) = \E(u) dx + \jmp{u}\odot \nu_u \H^{n-1}\restr J_u
\]
and in particular, $\Div u = \Tr\EE(u) = \Tr\E(u)dx + \jmp{u}\cdot\nu \H^{n-1}\restr J_u$.

   Since the singular part of the divergence of
$u$, $(\Div u)^s$, is given by $\ \jmp{u}\cdot \nu\H^{n-1}\restr J_u$,
the condition $\jmp{u}\cdot\nu\ge 0 \; \H^{n-1}-\textrm{a.e.\ in } J_u$  is equivalent to requiring that $\Div u^s$ be a nonnegative 
Radon measure, or, equivalently, that $\Dm u\in L^2(\Omega)$.

The goal of this paper is to show that, in the sense of $\Gamma$-convergence,
the energy $E(u)$ can be approximated with a sequence
of Ambrosio-Tortorelli{-type}~\cite{AT}  elliptic problems, given by
\begin{multline} \label{approxim}
E_\e(u,v) :=\\ {\frac12}\int_\Omega (\ee+v^2)
\left(2\mu |\EEd(u)|^2+ K(\Dp u)^2\right)
+  K(\Dm u)^2 dx  + G_c\int_\Omega \left\{\e|\nabla v|^2+\frac{(1-v)}{4\e}^2\right\}dx
\end{multline}
{where $\ee \ge 0$ is a parameter, see~\cite{AMM}.}

For simplicity 
%(and, in particular, the failed attempts 
%up to now to define correctly a space ``GSBD''),
we will assume that all functions satisfy a uniform,
{\it a priori} given $L^\infty$ bound. 
This is certainly a restriction because of the lack of a maximum principle in the context of linearized elasticity. We do not know at present how to remove this assumption. Given $M>0$, we introduce the following constrained functionals:
\[
E^M_\e(u,v)\ :=\begin{cases}
E_\e(u,v) & \ \textrm{ if } u\in H^1(\Omega;\R^n)\,, v\in H^1(\Omega),
\|u\|_{L^\infty}\le M\,, \\
+\infty & \ \textrm{otherwise,}
\end{cases}
\]
and:
\[
E^M_0 (u,v) \ :=\begin{cases}
E(u) &\ \textrm{ if } v=1\textrm{ a.e., }
u\in SBD(\Omega;\R^n), \|u\|_{L^\infty}\le M\,, \\
+\infty & \ \textrm{otherwise.}
\end{cases}
\]

%\end{remark}
Our result is as follows:
\begin{theorem}\label{th:main} Let the dimension {be} $n=2$.
Assume $\lim_{\e \to 0} \ee/\e=0$.
Then, $E^M_\e$ $\Gamma$-converges to $E^M_0$ as $\e\to 0$,
in $L^2(\Omega;\R^2)\times L^2(\Omega)$.
Moreover, if $(u_\e,v_\e)_{\e>0}$ is such that $\sup_{\e>0}E^M_\e{(u_\e,v_e)}<+\infty$,
then $\{(u_\e,v_\e):\e>0\}$ is sequentially precompact in $L^2$,
$v_\e\to 1$, and $E(u)\le \liminf_{k\to\infty} E_{\e_k}(u_{\e_k},v_{\e_k})$
for any limit point $u$ of a sequence $(u_{\e_k})$.
\end{theorem}
\begin{remark}
It will be clear from the proof that the result {also holds} if,
given $k>0$, $k\le K$, the terms $K(\Dp u)^2$  and $K(\Dm u)^2$
in~\eqref{approxim}
are replaced respectively with $(K-k)(\Div u)^2+ k(\Dp u)^2$ and
$k(\Dm u)^2$.
\end{remark}

\begin{remark}\label{rem:othermodels}
We  emphasize that, while  we present our
result and its proof in the case
of the simple constraint $\Div^s u\ge 0$,
the same proof  carries through for other constrained models such as
when only shear opening is present, as proposed in~\cite{LancioniRC}.  There the constraint
reads $\jmp{u}(x)\cdot\nu(x)= 0, \; \text{ for $\H^{n-1}$-a.e. } x\in J_u$.
{The opening constraint model for concrete~\cite{Ortiz-1985a,DelPiero-1989a}, implemented in~\cite{Freddi-Royer-Carfagni-2011a},
which boils down in the limit to $\jmp{u}(x)\in \R_+ \nu(x)$, is also {manageable}}. We refer to~\cite{FRC10} for a unified presentation of these cases.

In the first case, the approximating energy is
\begin{equation} \label{approxim-shear}
E_\e(u,v) = {\frac12}\int_\Omega 2\mu(\ee+v^2)
 |\EEd(u)|^2 + K (\Div u)^2 dx  + G_c\int_\Omega \left\{\e|\nabla v|^2+\frac{(1-v)}{4\e}^2\right\}dx
\end{equation}
and its $\Gamma$-limit -- still with the additional bound $\|u\|_\infty\le M$ --
will be the same as before {(see~\eqref{energy})}, but with the constraint
replaced with $\jmp{u}\cdot \nu(x)=0$ a.e.~on $J_u$.
In the second case, {a possible approximation is}
\begin{multline} \label{approxim-positive}
E_\e(u,v) = {\frac12}\int_\Omega (\ee+v^2)\Big(2\mu
 |\EE^+(u)|^2 +\lambda (\textup{Tr\,}\EE^+(u))^2\Big)
 + 2\mu|\EE^-(u)|^2 +\lambda (\textup{Tr\,}\EE^-(u))^2dx 
\\ + G_c\int_\Omega \left\{\e|\nabla v|^2+\frac{(1-v)}{4\e}^2\right\}dx
\end{multline}
where $\EE^+(u)$ is the projection of $\EE(u)$ {onto} the cone of nonnegative symmetric matrices
and $\EE^-(u)=\EE(u)-\EE^+(u)$.
The $\Gamma$-limit of this energy
is given by the same functional $E(u)$.  {The constraint is now}
that the singular part $\EE^s(u)$ should only have non-negative eigenvalues. {Since}
$\EE^s(u)=\jmp{u}\odot \nu\H^{n-1}\restr J_u$ is rank-1-symmetric
($\jmp{u}\odot \nu=(\jmp{u}\otimes \nu + \nu\otimes\jmp{u})/2$,
hence has rank 1 or 2), this implies
that $\jmp{u}$ and $\nu$ should be a.e.~aligned and in the same direction.
{Indeed, it is easy to see that given two vectors $a,b\in \R^2$, then
$\textup{det} (a\otimes b+b\otimes a) = -(a\times b)^2/4\le 0$ so that
this matrix is nonnegative only if the vectors are aligned and in
the same direction.}
\end{remark}

\begin{figure}[htb]
\begin{center}
\hspace{6mm}\raisebox{8mm}{\includegraphics[width=.4\textwidth,trim={0 0 22.5cm 18.8cm},clip]{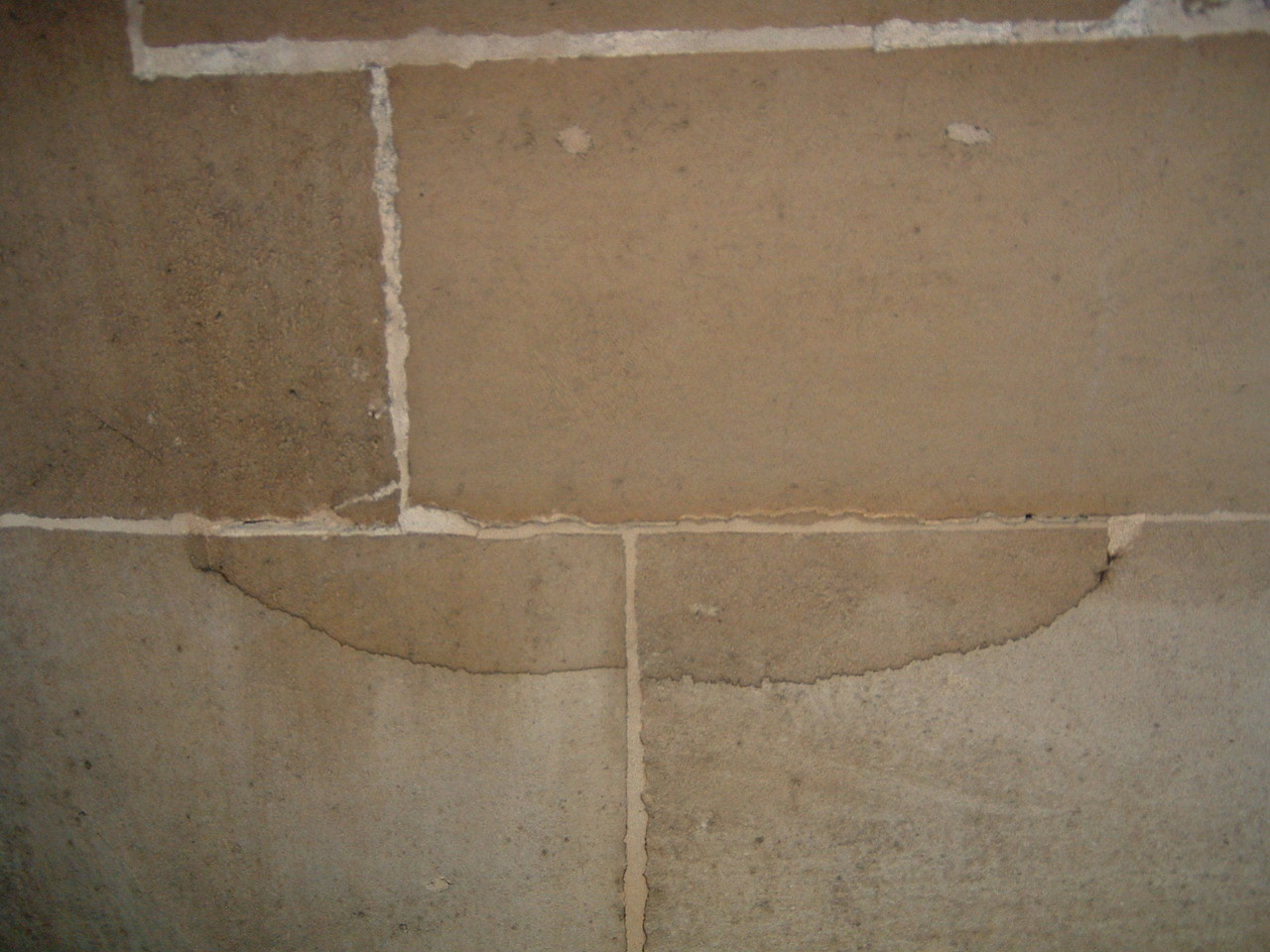}}
\hfil
\includegraphics[width=.5\textwidth]{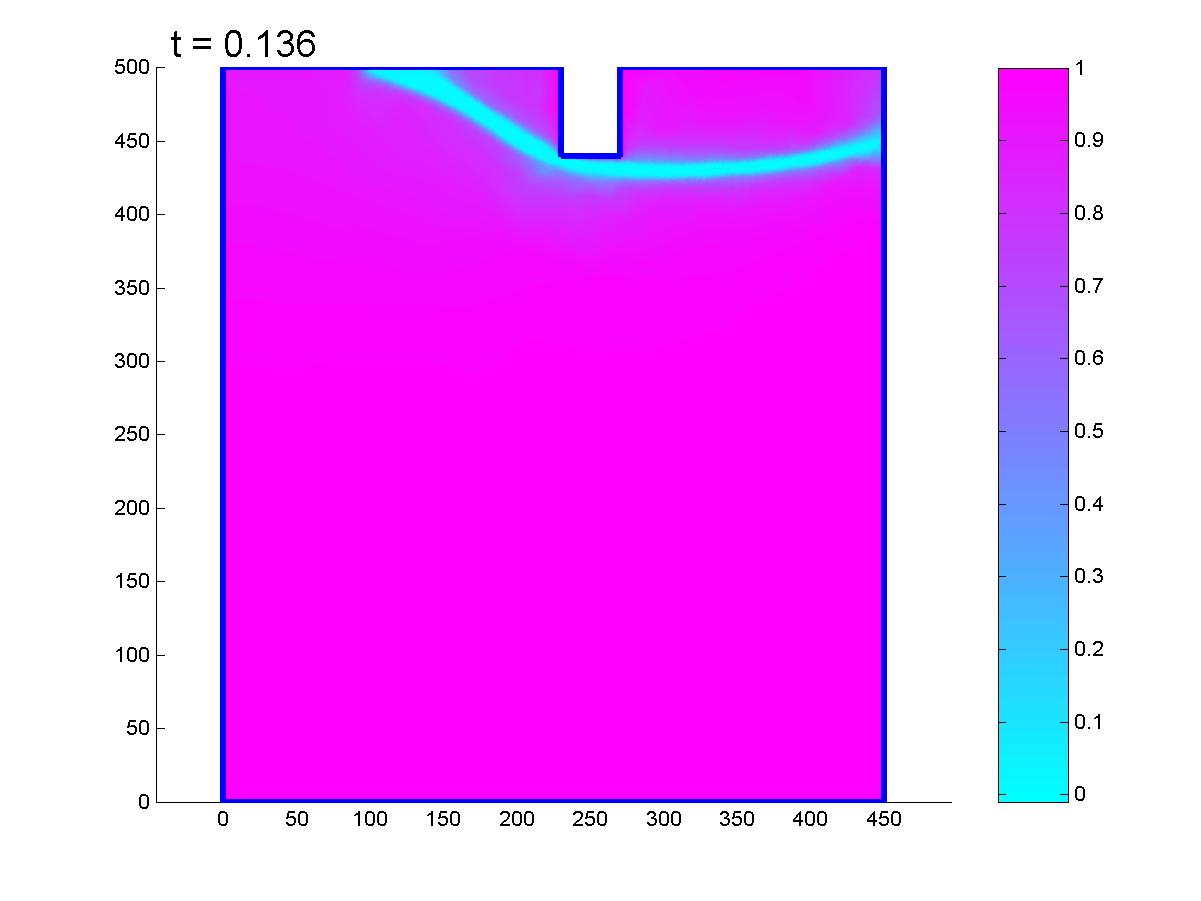}
\end{center}
\caption{On the left, a shear fracture observed on
a stone of the Pantheon in Paris. Right, a simulation based on the
energy~\eqref{approxim-shear}, computed by the authors of~\cite{LancioniRC}.}\label{fig:LRC}
\end{figure}

Although the mathematical proof
of the relevance of such an approximation was wanting up till now,
the numerical pertinence  of the approximating functional given in Theorem~\ref{th:main} \blue{or its variants} for dealing with non-interpenetration  has  been successfully demonstrated in e.g.~\cite{AMM,Freddi-Royer-Carfagni-2011a,LancioniRC}; see Figure~\ref{fig:LRC} for  illustration.

\blue{
\begin{remark}
The techniques developed over the years for this paper (which is part of
a project initiated more than 10 years ago),
starting from the Korn-Poincar\'e inequality in~\cite{CCF_KPI},
have subsequently been successfully adapted to the study of other
problems involving $GSBD$ functions,
such as a new approximation result~\cite{ChambolleCrismale2017},
or the proof of existence of strong minimizers for Griffith type
energies in higher dimensions~\cite{ChambolleContiIurlano}.
\end{remark}
}

The organization of the paper is as follows. %[TO BE DESCRIBED HERE]
We first give a proof of our main result {along the  lines of a classical argument in~\cite{AT} which
have since been reproduced and adapted to many settings}.
However, the proof of the $\Gamma$-limsup, which we sketch, relies on {an  approximation of the limit function $u$ with a function
exhibiting a ``simple jump''. Such a result is lacking at present.
We thus follow a different strategy described} in
Section~\ref{secfullproof}. Most of {the} proof can be carried {out} in arbitrary dimension, {although,}
 as explained at the very end of Section~\ref{sec:GenProof}, we lack a crucial
estimate to conclude. The last section shows how we can circumvent this difficulty in 
dimension $n=2$.

In the rest of the paper, we will assume, 
{in the sake of simplicity}, that $\mu=1$, {$K=2$}, $G_c=1${; the proof
 clearly does not depend} on the values of these  parameters.

\section{A first, partial proof of convergence}\label{secpartialproof}

In order to establish Theorem~\ref{th:main}, we need to show that
\begin{itemize}
\item[(i)] If $\sup_{k\ge 1} E^M_{\e_k}(u_k,v_k)<+\infty$
for a given sequence $\e_k\downarrow 0$, then $v_k \to 1$,
$u_k$ converges in $L^2(\Omega;\R^n)$ to some displacement $u$ (up
to a subsequence, and with obviously $\|u\|_{L^\infty}\le M$),
and  $$E(u)\le \liminf_{k\to\infty}E^M_{\e_k}(u_k,v_k)$$
(compactness and $\Gamma$-lim\,inf inequality);
\item[(ii)] For each $u\in SBD(\Omega)$, with $\|u\|_{L^\infty}\le M$\,,
there exists $(u_\e,v_\e)$ with $\lim_{\e\to 0}\|u_\e-u\|_{L^2}=0$
and $$\limsup_{\e\to 0} E^M_\e(u_\e,v_\e)\le E(u)$$ ($\Gamma$-lim\,sup
inequality).
\end{itemize}

In the next Subsection, we quickly establish the first point   mostly
following and detailing the proof in~\cite{ChaSBD1}.
Then, in Subsection~\ref{secfirstlimsup}, we propose a proof of~(ii),
still inspired from~\cite{ChaSBD1}, but valid only for a subclass of
$SBD$ displacements. The complete proof is given in Section~\ref{secfullproof}.
{Since large parts of the argument hold in any dimension, we formulate them for a general $n\ge2$, and only use the assumption $n=2$ in  the final construction discussed in Section \ref{sec:2-d-case}.}

\subsection{Compactness and proof of the liminf inequality in  (i)}\label{secliminf}

The proof of the liminf inequality is {that of} the standard
case, {since} the new functional is larger while its limit, at least
on its domain of definition, is the same. {The proof  {detailed} below
is} adapted from~\cite{ChaSBD1}.

Let $\e_k>0$ be a sequence converging monotonically to zero, and
let $u_k,v_k$ be a displacement and a function such that
$E_{\e_k}(u_k,v_k)\le C<+\infty$. First, we observe
that since
\[
\int_\Omega (1-v_k)^2\,dx \ \le\ \blue{4C\e_k}\,,
\]
$v_k\to 1$ in $L^2(\Omega)$. We also have
\[
\int_\Omega \e_k|\nabla v_k|^2\,+\,\frac{1}{4\e_k}
(1-v_k)^2\,dx\ \ge\ \int_\Omega |1-v_k||\nabla v_k|\,dx\,,
\]
so that, using the co-area formula, %and increasing rearrangement,
 we find that
\begin{multline}\label{Eksplit}
E_{\e_k}(u_k,v_k) \ge \\ \int_0^1
\int_{\{v_k>s\}} \hspace{-4mm}\left( 2s(|\EEd(u_k)|^2+(\Dp u_k)^2) + (\Dm u_k)^2 dx
+ (1-s)\H^{n-1}(\partial^* \{v_k>s\})\right) ds
\end{multline}
where $\partial^* \{v_k>s\}$ denotes the reduced boundary of the superlevel set 
$ \{v_k>s\}$.

First, we deduce that for each $k$, we can find $s_k\in (1/4,3/4)$
such that the function $\tilde{u}_k{:=}u_k\,\car{\{v_k>s_k\}}$ is in
$SBD(\Omega)$, with $\|\tilde{u}_k\|_\infty\le M$, $\blue{J_{\tilde u_k}}=\partial^*\{v_k>s_k\}$
and
\[
\frac{1}{2} \int_\Omega |\E(\tilde{u}_k)|^2\, +\,
\frac{1}{4}\H^{n-1}(J_{\tilde{u}_k}) \ \le \ C \,<\,+\infty\,.
\]
>From the compactness Theorem in~\cite[Thm~1.1]{BCDM}, we deduce that
up to a subsequence, $\tilde{u}_k$ {converges}, in $L^2$,
to some $u\in SBD(\Omega)$,
with $\E(u)\in L^2$ and $\H^{n-1}(J_u)<+\infty$. Now, since
the sequence $(u_k)$ is uniformly bounded in the $L^\infty$ norm,
and $v_k\to 1$ in $L^2(\Omega)$, so that $|\{v_k\le s_k\}|\to 0$,
we deduce that $u_k\to u$ in $L^2(\Omega)$ (or any $L^p$, for $p<+\infty$).

Now, for a.e. $s\in (0,1)$, one must have
$u_k\car{\{v_k>s\}}\to u$, in $L^2$, as $k\to\infty$.
% From this, we get
The $SBD$ variant of Ambrosio's compactness and semicontinuity
theorem, proved in~\cite{BCDM}, yields that
\begin{multline*}
\liminf_{k\to\infty}\int_{\{v_k > s\}}
\left( 2s(|\EEd(u_k)|^2+(\Dp u_k)^2\right)  + (\Dm u_k)^2\,dx
+ (1-s)\H^{n-1}(\partial^*\{v_k>s\})
\\ \ge
\int_\Omega 2s(|\Ed(u)|^2+(\Trp\E(u))^2) + (\Trm\E(u))^2\,dx + (1-s)\H^{n-1}(J_u).
\end{multline*}
 { Here $\Trp$ and $\Trm$ denote} the positive part and
the negative part of the trace, respectively.

However, integrating back this inequality with respect to $s\in [0,1]$
does not allow {one} to recover~(i). Indeed, the weight in front
of the surface term will only be $\int_0^1 (1-s)ds=1/2$. This is
because in the semicontinuity result, one loses the fact that
the jump set is obtained as the limit of the collapsing level sets
$\{v_k\le s\}$, and therefore the perimeter of these {sets} \blue{measures}
twice the size of the limiting jump set.
This heuristic observation is easy to actually turn into
a proof. It follows a variant of Ambrosio's theorem which is found 
in~\cite[Thm~2]{BCS} (see also Lemma~2 in~\cite{BC}).
It is written there for scalar or vectorial $GSBV$ functions, but its proof,
which is based on slicing, can easily be shown to extend to similar problems in $SBD$
(with an $L^\infty$ bound). It shows that in fact, for a.e.~$s$,
\begin{multline*}
\liminf_{k\to\infty}\int_{\{v_k > s\}}
\left( 2s(|\EEd(u_k)|^2+(\Dp u_k)^2\right)  + (\Dm u_k)^2\,dx
+ (1-s)\H^{n-1}(\partial^*\{v_k>s\})
\\ \ge
\int_\Omega 2s(|\Ed(u)|^2+(\Trp\E(u))^2) + (\Trm\E(u))^2\,dx + 2(1-s)\H^{n-1}(J_u).
\end{multline*}
Using~\eqref{Eksplit} and Fatou's lemma, we deduce that
\begin{multline*}
\liminf_{k\to\infty} E_{\e_k}(u_k,v_k)\\
\ge \int_0^1 \left(
\int_\Omega 2s(|\Ed(u)|^2+(\Trp\E(u))^2) + (\Trm\E(u))^2\,dx + 2(1-s)\H^{n-1}(J_u)
\right)
\\
=\int_\Omega |\Ed(u)|^2+(\Tr\E(u))^2 dx+ \H^{n-1}(J_u).
\end{multline*}
{Finally observe that since $\Div u_k\wto \Div u$ as measures
and $\Dm u_k$ is bounded in $L^2(\Omega)$, denoting $f\in L^2(\Omega)$
a $L^2$-weak limit point of $\Dm u_k$, we deduce from 
the inequality $-\Div u_k\le \Dm u_k$ that $-\Div u\le f$, }
showing that $\Dm u\in L^2(\Omega)$ and therefore that $\jmp{u}\cdot \nu\ge 0$
$\H^{n-1}$-a.e.~on $J_u$.
The proof of compactness and the lower $\Gamma$-limit estimate
is complete.

\subsection{{A first proof of the limsup inequality in (ii),
when the jump set is ``nice'' enough}}
\label{secfirstlimsup}
A ``standard'' proof of a result such as Theorem~\ref{th:main}
would now show inequality (ii) first for ``simple'' displacements
(for instance, with smooth jump sets), and then show that doing so is not restrictive
by constructing, for an arbitrary $SBD$ displacement $u$, a sequence of
approximate ``simple" displacements $u_n$ with $E(u_n)\to E(u)$.
A diagonalization process would then be invoked to deduce (ii) in the general case.
This is for instance what is done in~\cite{ChaSBD1},
where such an approximation
is provided.
However, that particular approximation does not enjoy the constraint
$\jmp{u}\cdot\nu\ge 0$
on the jump set, and it is far from clear how to modify it to ensure this
constraint.

The {\it bona fide} proof of estimate (ii) is quite involved;
see Section~\ref{secfullproof}.
Below we provide an elementary (and
classical) proof in the particular case where $J_u$ 
is essentially closed, {\it i.e.}, $\H^{n-1}(\bar J_u\setminus J_u)=0$ and satisfies a lower density bound 
\begin{equation}\label{ldb}\H^{n-1}(J_u\cap B(x,r))\ge \kappa r^{n-1}, \;x\in \Omega, r< \text{dist }(x,\partial\Omega),
\end{equation}
so that  its
Hausdorff measure is given by its Minkowski content
\begin{equation}\label{eqMink}
\lim_{t\to 0} \frac{\left\{x\in\Omega\,:\, \dist(x,J_u)\le t\right\}}{2t}
=\H^{n-1}(J_u)
\end{equation}
(see for instance~\cite[Subsection 2.13]{AFP}).  
It is well known that $SBD_2$ fields can be approximated in energy
by fields which satisfy these conditions~\cite{ChaSBD1,ChaSBD2,Iurlano}. {However, none of the known
constructions ensure that a  constraint such as $\jmp{u}\cdot\nu\ge 0$ can be}
maintained in the approximation.

We also assume that $M=\|u\|_{L^\infty}<+\infty$. We then
choose $\phi$ a symmetric mollifier with support in $B(0,1)$. We
let $\delta_\e=\sqrt{\e\eta_\e}$ be an intermediate scale between
$\eta_\e$ and $\e$,
set as usual $\phi_{\delta_\e}(x):=(\delta_\e)^{-n}\phi(x/\delta_\e)$,
and define
\begin{equation*}
v_\e(x):= \gamma\left(
\frac{(\dist(x,J_u)-\delta_\e)^+}{\e}\right),
\quad 
u_\e :=\phi_{\delta_\e}*u
\end{equation*}
with  $u$  extended slightly out of $\Omega$, as explained farther at the onset of Subsection~\ref{sec:GenProof}.
Here $\gamma: [0,+\infty)\to [0,1]$ is the one-dimensional optimal
profile associated to the energy $\int_\Omega\e |\nabla v|^2+(1-v)^2/(4\e)dx$,
that is $$\gamma(t)=1-\exp(-t/2).$$

Then, $v_\e\to v=1$ and $u_\e\to u$ in $L^2$. On the other hand
\[
\EE(u_\e)(x) = \phi_{\delta_\e}*\E(u)(x)
+\int_{J_u} \phi_{\delta_\e}(x-y)\jmp{u}(y)\odot \nu_u(y)d\H^{n-1}(y)
\]
so that if $\dist(x,J_u)>\delta_\e$, $\EE(u_\e)=\phi_{\delta_\e}*\E(u)$
while in general, $|\EE(u_\e)|\le cM/\delta_\e$ for some constant
$c>0$ depending only on $\phi$ (and $n$). Hence,
\begin{multline}\label{estimbulk}
\int_\Omega (\eta_\e + v_\e^2) (|\EEd(u_\e)|^2 + (\Dp u_\e)^2) + (\Dm u_\e)^2dx
\\
\le (1+\eta_\e) \int_{\{\dist(\cdot,J_u)> \delta_\e\}} |\phi_{\delta_\e}*\Ed(u)|^2
+(\phi_{\delta_\e}*\Div u)^2 dx
\\
+\int_{\{\dist(\cdot,J_u)\le \delta_\e\}} (\Dm u_\e)^2 dx+c|\{\dist (\cdot,J_u)\le\delta_\e\}|\eta_\e \frac{M^2}{\delta_\e^2}.
\end{multline}
Since 
$$
0\le \Dm u_\e = (\phi_{\delta_\e}*\Div u)^-\le {\phi_{\delta_\e}} * \Dm u
$$
 and because the latter
is uniformly bounded in $L^2$, and since further
$|\{\dist(\cdot,J_u)\le \delta_\e\}|\to 0$ we deduce that
\[
\limsup_{\e \to 0} \int_{\{\dist(\cdot,J_u)\le \delta_\e\}} (\Dm u_\e)^2dx=0.
\]
Therefore, recalling that, thanks to Eq.~\eqref{eqMink},
$|\{\dist(\cdot,J_u)\le \delta_\e\}|\approx 2\delta_\e\H^{n-1}(J_u)$ as
$\e\to 0$ and $\eta_\e/\delta_\e\to 0$,
Eq.~\eqref{estimbulk} becomes in the limit
\begin{equation}\label{eq:lsu}
\limsup_{\e\to 0}
\int_\Omega (\eta_\e + v_\e^2) (|\EEd(u_\e)|^2 + (\Dp u_\e)^2) + (\Dm u_\e)^2dx
\le\int_\Omega (|\Ed(u)|^2+(\Tr\E(u))^2)\,dx.
\end{equation}

On the other hand, since $|\nabla \dist(\cdot,J_u)|=1$ a.e.,
\[
\e |\nabla v_\e|^2=\frac{(1-v_\e)^2}{4\e}=\frac{1}{4\e}
\exp\left(-\frac{(\dist(x,J_u)-\delta_\e)^+}{\e}\right)
\]
a.e.~in $\{\dist(\cdot,J_u)\ge \delta_\e\}$ hence, using the co-area formula,
\begin{multline*}
\int_\Omega \e|\nabla v_\e|^2+\frac{(1-v_\e)}{4\e}^2 dx=
\frac{|\{\dist(\cdot,J_u)< \delta_\e\}|}{4\e} +
\frac{1}{2\e}\int_{\{\dist(\cdot,J_u)\ge \delta_\e\}} e^{-\frac{(\dist(x,J_u)-\delta_\e)^+}{\e}}dx
\\
= \frac{|\{\dist(\cdot,J_u)< \delta_\e\}|}{4\e} +
\frac{1}{2\e}\int_{\delta_\e}^\infty
e^{-\frac{s-\delta_\e}{\e}}\H^{n-1}(\partial \{\dist(\cdot,J_u){<s}\}\cap\Omega)ds.
\end{multline*}
Let $f(s):=|\{\dist(\cdot,J_u)<s\}|$. By the co-area formula,
$f'(s)=\H^{n-1}(\partial \{\dist(\cdot,J_u){<s}\}\cap\Omega)$ for a.e.~$s>0$.
We find that
\begin{multline*}
\int_\Omega \e|\nabla v_\e|^2+\frac{(1-v_\e)}{4\e}^2 dx= \frac{f(\delta_\e)}{4\e}
+\frac{1}{2\e}\int_{\delta_\e}^\infty e^{-\frac{s-\delta_\e}{\e}} f'(s)ds \\
= -\frac{f(\delta_\e)}{4\e}+
\frac{1}{2\e} \int_{\delta_\e}^\infty e^{-\frac{s-\delta_\e}{\e}} f(s) ds
\\
= -\frac{\delta_\e}{\e}\frac{f(\delta_\e)}{4\delta_\e}+e^{\frac{\delta_\e}{\e}}\int_{\delta_\e/\e}^\infty s e^{-s} \frac{f(\e s)}{2\e s}ds.
\end{multline*}
By~\eqref{eqMink}, $f(\e s)/(2\e s)\to \H^{n-1}(J_u)$ as $\e\to 0$ for all $s>0$.
In general, denoting $\ell=\lim_{\e\to 0}\delta_\e/\e$
and assuming {$\ell<+\infty$}, we find in the limit
\begin{multline}\label{eq:lsv}
\limsup_{\e\to 0} \int_\Omega \e|\nabla v_\e|^2+\frac{(1-v_\e)}{4\e}^2 dx
\\
\le
\H^{n-1}(J_u)\left(-\frac{\ell}{2}+e^\ell\int_\ell^\infty s e^{-s}\,ds\right)
=\left(1+\frac{\ell}{2}\right)\H^{n-1}(J_u).
\end{multline}
Since we have assumed {$\delta_\e=\sqrt{\eta_\e\e}$ we have} $\ell=0$ and the right-hand side is
simply $\H^{n-1}(J_u)$.
Collecting~\eqref{eq:lsu} and~\eqref{eq:lsv} yields the desired estimate.

\section{{A general proof in dimension 2}}\label{secfullproof}

We now describe the general proof of the lim-sup inequality, which
will work without further {assumptions} on the jump set of $u$. 
Most of the proof can be carried {out} in arbitrary dimension $n\ge 2$.
A technical difficulty {will prevent us from  concluding when $n>2$. }

We choose a small parameter $\theta\ll1$. 
Since the jump set $J_u$ is countably rectifable, there exists a regular part $\Gamma$, a finite union of closed subsets of  $\mathcal C^1$-hypersurfaces  such that $\H^{n-1}(J_u\triangle \Gamma)\le \theta^2$. 

 We fix a small length scale $\delta:=\e\ell$, {$\ell\in(0,1)$} small,
and subdivide the domain into cubes {$Q_z$} of size proportional to $\delta$ (details below). 
We call  {$Q_z$} {\it {\it good}  } if it  contains an amount of jump smaller than $\theta\delta^{n-1}$, that is if 
\begin{equation}\label{goodcube}\H^{n-1}(J_u\cap Q_z)\le\theta \delta^{n-1};\end{equation}
 otherwise  {$Q_z$} is {\it bad}.

The function $v$ is constructed so that it vanishes on  a $\delta$-neighbourhood $\Sigma$ of both $\Gamma$ and  the {\it bad} cubes.

In the {\it bad} cubes we shall use a mollification of $u$, in the {\it good}   ones a mollification of $u$ after ``cleaning out'' the small holes
using the rigidity {result} of~\cite[Prop.~3.1]{CCF_KPI}. The result reads as follows:

\begin{proposition}\label{th:prop3b} 
Let %$0<\theta''<\theta'<1$,
 $r>0$. 
Let $Q=(-r,r)^n$,
	$Q'=(-r/2,r/2)^n$, $p\in [1,\infty)$, $u\in SBD^p(Q)$.
	\begin{enumerate}
		\item \label{propsobol1}
		There exists a set $\omega\subset Q'$ and an affine
		function $a:\R^n\to\R^n$ with $\EE(a)=0$ such that
		\begin{equation}\label{eq:stimom}
		|\omega|\le c_* r \H^{n-1}(J_u)
		\end{equation}
		and
		\begin{equation}\label{eq:stimu}
		\int_{Q'\setminus\omega} |u-a|^{np/(n-1)}{dx}
		\le c_*r^{n(p-1)/(n-1)}\left(\int_Q |\E(u)|^p dx\right)^{n/(n-1)}.
		\end{equation}
		\item \label{propsobol2}
		If additionally $p>1$ then there is $\bar p>0$ (depending on $p$ and $n$) such that,
		for a given mollifier $\phi\in C^\infty_c(B_{1/2})$ with
                $\int \phi\,{dx}=1$, letting  $\phi_r(x)=r^{-n}\phi_1(x/r)$,
		the function $v=u\chi_{Q'\setminus\omega}+a\chi_\omega$ obeys
		\begin{equation*}
		\int_{Q''} | \EE(v\ast \phi_r)- \E(u)\ast\phi_r|^p dx \le c 
		\left(\frac{\H^{n-1}(J_u)}{r^{n-1}}\right)^{\bar p} \int_Q |\E(u)|^p dx\,,
		\end{equation*}
		where $Q''=(-r/4,r/4)^n$.
	\end{enumerate}
	The constant in \blue{(\ref{propsobol1})} depends only on $p$, $n$ 
        the one in \blue{(\ref{propsobol2})} also on $\rho$. % and $\theta''$.
\end{proposition}
\begin{remark}\label{rem:linfty}
Thanks to Lemma~\ref{lemmalinfty} in the Appendix, one can assume
additionally that $\|a\|_{L^\infty(Q')}\le \|u\|_{L^\infty(Q)}$ if in addition,
$u$ is bounded.
\end{remark}

 The challenge {-- which unfortunately
we cannot overcome  except in 2D --}
 will be in the handling of   the boundary between  {\it the good}   and {\it the bad} regions.
Hence in a second step, we shall further introduce ``boundary {\it good}   cubes''.
On these we can clean up the jump before mollification
{ using a construction due to~\cite{CFI}, which we only know to hold true
in dimension 2. The details are found in Subsection \ref{sec:2-d-case}.}

\subsection{{The general proof}}
\label{sec:GenProof}
We first assume that $u$ is defined slightly outside of $\Omega$
in a domain $\Omega'\supset\supset\Omega$. The necessary assumption
for this is that
$\partial\Omega$ be  a subgraph locally: then,
the construction consists in translating
 $u$  outside of $\Omega$ near the boundary and in 
glueing
the pieces together with a partition of unity. This creates a new $u'$ with still,
$\jmp{u'}\cdot\nu \ge 0$ (or $=0$) on $J_{u'}$
{and $\|u'\|_{L^\infty(\Omega')}\le \|u\|_{L^\infty(\Omega)}\le M$.}
Also, we can assume
$\H^{n-1}(\partial \Omega\cap J_{u'})=0$. We drop the ``prime''
and denote the extended function by $u$ in the following. As usual we also set
$\Omega^\delta :=\{x\in\R^n\,:\,\dist(x,\Omega)<\delta\}\subset\subset\Omega'$
for $\delta>0$ small enough.

{We} consider the cubes $Q_z = z+(-4\delta,4\delta)^n$, $\tilde q_z = z+(-2\delta,2\delta)^n$, $q_z = z+(-\delta,\delta)^n$, for $z\in (2\delta)\Z^d$.
We also consider $\phi_\delta(x) = \delta^{-n}\phi(x/\delta)$ a mollifier
with support in $B(0,\delta/2)$.

We let
\begin{align*} \Omega_g^\delta := \bigcup\{ q_z\,:\, Q_z \subset\Omega'\textup{ \it\; good } \},\\
\Omega_b^\delta := \bigcup\{ q_z\,:\, Q_z \subset\Omega'\textup{ \it \; bad} \}.
\end{align*}

We set (almost) as before, for $x\in\Omega'$,
\[
v^0_\e(x):= \gamma\left(
\frac{(\dist(x,\Gamma)-16\sqrt{n}\delta)^+}{\e}\right)
\]
and find as before (\textit{cf} Eq.~\eqref{eq:lsv}), {recalling that $\delta=\e\ell\le\e$},
that
\begin{equation}\label{surfest-main}
\limsup_{\e\to 0} \int_\Omega \e|\nabla v^0_\e|^2+\frac{(1-v^0_\e)}{4\e}^2 dx
\le
\left(1+8\sqrt{n}\ell\right)\H^{n-1}(\Gamma\cap \Omega).
\end{equation}

\begin{remark} If $Q_z$ intersects $\Gamma$ then all points in $Q_z$ are at a distance less than {$8\sqrt n \delta$} of $\Gamma$ so that  $v^0_\e\equiv 0$ on $Q_z$.
\end{remark}

We wish to define $v_\e$ as   zero only near $\Gamma$ and around the {\it bad} cubes. For those {\it bad} cubes that intersect $\Gamma$ we  take $v_\e^z:= v^0_\e$. For those that do not intersect $\Gamma$, we take 
 $$v^z_\e(x)= \left\{\begin{array}{ll}0 &\mbox{ in }
 B(z,16\sqrt{n}\delta)\supset Q_z\\ 
 \gamma((|x-z|-16\sqrt{n}\delta)^+/\e) &\mbox{ else}.
\end{array}\right.$$

A simple calculation would show the existence of a constant $C$ such that
\begin{equation}\label{surfest}
\int_\Omega \e|\nabla v^z_\e|^2+\frac{(1-v^z_\e)}{4\e}^2 dx
\le C\e^{n-1}.
\end{equation}

Denoting {by} $B\!C$ the set of the {\it bad\/} cubes  that do not intersect $\Gamma$, its cardinality % (denoted by $\#B\!C$)
satisfies
\begin{equation}\label{badcubes}
\#B\!C\le C\H^{n-1}(J_u\setminus \Gamma)/ (\theta\delta^{n-1})\le C\theta/\delta^{n-1}.
\end{equation}
Indeed, in view of Eq.~\eqref{goodcube}, 
$$
\theta\delta^{n-1} \#B\!C\le \sum_{\text{{B\!C}}} \H^{n-1}(J_u\cap Q_z)=\sum_{\text{{B\!C}}} \H^{n-1}(J_u\cap Q_z\setminus \Gamma).
$$
But there is at most $C$ (some constant) overlaps  between those cubes so that, in view of the above,
$$
\theta\delta^{n-1} \#B\!C\le C \H^{n-1}(J_u\setminus \Gamma)\le C \theta^2.
$$

In view of eqs~\eqref{surfest-main},\eqref{surfest},\eqref{badcubes}, if $v_\e^\ell$ is the
{min} of $v^0_\e$ and all the $v^z_\e$
defined for the {\it bad} cubes not intersecting $\Gamma$, then
\begin{equation*}
\limsup_{\e\to 0} \int_\Omega \e|\nabla v_\e^\ell|^2+\frac{(1-v_\e^\ell)}{4\e}^2 dx
\le
\left(1+8\sqrt{n}\ell\right)(\H^{n-1}(J_u\cap \Omega)+\theta^2)
+ C\theta \frac{1}{\ell^{n-1}},
\end{equation*}
hence
\begin{equation}\label{eq:estimvpart}
\limsup_{\ell\to 0}\limsup_{\theta\to 0}\limsup_{\e\to 0} \int_\Omega \e|\nabla v_\e^\ell|^2+\frac{(1-v_\e^\ell)}{4\e}^2 dx
\le
\H^{n-1}(J_u\cap \Omega).
\end{equation}
This takes care of the surface term.

{We remark in addition that, by construction, {
\begin{equation}\label{eq:vol-v=0}
 |\{v_\e^\ell=0\}|=O(\delta).
 \end{equation} 
 Specifically that set} can be decomposed as the union of
$\{v^0_\e=0\}$, which has a volume of order $C\delta\H^{n-1}(\Gamma)$, and
of the union of $\{v^z_\e=0\}$ for all the bad cubes in $B\!C$, which
has a volume of order $C\theta\delta$.}
%and $v_\e$ is nonzero only on {\it good}   cubes $q_z$ whose neighbouring $q_{z'}$ are also {\it good}  .

\vskip1cm

Take any of the cubes $Q_z$. 
{From Proposition~\ref{th:prop3b}
 \cite[Prop~3.1]{CCF_KPI}},
in $Q_z$ there exists $a_z:\R^n\to \R^n$ affine with $\EE(a_z)=0$ and $\omega_z\subset \tilde q_z$ such that 
$|\omega_z|\le \blue{c_*}\delta\H^{n-1}(J_u\cap Q_z)$, and % for $p=2n/(n-1)>2$,
\begin{equation}\label{eq:estimaz}
\left(\int_{\tilde q_z\setminus \omega_z} |u-a_z|^{\frac{2n}{n-1}}\,dx\right)^{1-1/n}
 \le \blue{c_*} \delta\int_{Q_z} |\E(u)|^2\,dx,
\end{equation}
while, moreover, setting 
\begin{equation}\label{eqdefwz}
w_z := u\chi_{\tilde q_z\setminus \omega_z}+ a_z\chi_{\omega_z} 
\end{equation}
one has,
for some $q=q(n)>0$ {and for a given mollifier $\phi\in C^\infty_c(B_{1/2})$ with
                $\int \phi{\,dx}=1$},
\begin{equation}\label{estimdif}
\int_{q_z} |\EE(w_z*\phi_\delta)-\E(u)*\phi_\delta|^2 \,dx
\le c\left(\frac{\H^{n-1}(J_u\cap Q_z)}{\delta^{n-1}}\right)^q
\hspace{-.5em}\int_{Q_z} |\E(u)|^2\,dx.
\end{equation}
Finally, 
thanks to Remark~\ref{rem:linfty}, it is also possible to assume {that}
\begin{equation}\label{linfbound}
\|a_z\|_{L^\infty(\omega_z)}\le\|a_z\|_{L^\infty(\tilde q_z)}\le \|u\|_{L^\infty(\Omega)}.
\end{equation}

We first work with the {\it good} cubes,
 and more precisely we restrict this terminology to the {\it good} cubes $Q_z$ { on which}
 $v^\ell_\e\not\equiv 0$ (hence not too close to $\Gamma$ or a bad cube), modifying accordingly the definition of the sets $\Omega_g^\delta$, $\Omega_b^\delta$.
 An observation is that if $Q_z,Q_{z'}$ are  two  {\it good}   cubes 
{such that $q_z,q_{z'}$ are touching}
(by which we mean  that $|z-z'|_\infty = {2\delta}$), then the volume of $\tilde q_z\cap \tilde q_{z'}$
is at least $\delta^n$. Furthermore, since $a_z$ and $a_z'$ are affine,
$$
|(\tilde q_z\setminus \omega_z)\cap (\tilde q_{z'}\setminus \omega_{z'})|{ \|a_z-a_{z'}\|_{L^\infty(Q_z\cap Q_{z'})}^{\frac{2n}{n-1}}}\le C\int_{(\tilde q_z\setminus \omega_z)\cap (\tilde q_{z'}\setminus \omega_{z'})}|a_z-a_z'|^{\frac{2n}{n-1}}\,dx.
$$
It then follows from~\eqref{eq:estimaz} that, for some constant $C$,
\[
|(\tilde q_z\setminus \omega_z)\cap (\tilde q_{z'}\setminus \omega_{z'})|^{1-1/n}{ \|a_z-a_{z'}\|_{L^\infty(Q_z\cap Q_{z'})}^2}  \le 
C\delta \int_{Q_z\cup Q_{z'}}|\E(u)|^2dx
\]
so that, because $Q_z, Q_z'$ are {\it good} cubes,
\begin{equation}\label{eq:diffazz}
{
\|a_z-a_{z'}\|_{L^\infty(Q_z\cap Q_{z'})}^2} \le \frac{C \delta}{(\delta^n(1-c\theta))^{\frac{n-1}{n}}}\int_{Q_z\cup Q_{z'}}|\E(u)|^2dx
\le \frac{C}{\delta^{n-2}}\int_{Q_z\cup Q_{z'}}|\E(u)|^2dx.
\end{equation}
if $\theta$ is small enough.
We can order (arbitrarily) all $z$ such that $Q_z$ is {\it good}   and
denote the corresponding sequence $\{z_j\}_{j\in \GC}$, where $\GC$ denotes a numeration of the {\it good} cubes. Then, we define
\[
\tilde u (x) = \begin{cases}
u(x) & \textup{ if } x\in   
\Omega' \setminus \bigcup_{j\in \GC} \omega_{z_j} \\
a_{z_j} & \textup{ if } x \in \omega_{z_j} \setminus \bigcup_{i<j} \omega_{z_i}.
\end{cases}
\]
Observe that thanks to~\eqref{linfbound}, $\|\tilde u\|_{L^\infty(\Omega')}\le \|u\|_{L^\infty(\Omega)}\le M$.
We let $u_\e := \tilde u*\phi_\delta$, and, in order to provide an estimate for  the volume term,
we now propose to bound from above $\int_{\Omega_g^\delta}|\EE(\tilde u_\e)|^2\,dx$ . This is done
by showing that this function is, in $L^2(\Omega_g^\delta)$, close to $\E(u)$.

Upon decomposing $\tilde u$ as $\tilde u=w_{z_j}+(\tilde u-w_{z_j})$, \blue{with $w_{z_i}$ defined in (\ref{eqdefwz}),} we obtain:
\begin{multline}\label{eq:globalqzdiff}
\int_{\Omega_g^\delta}|\EE(u_\e)-\phi_\delta*\E(u)|^2\,dx
\\\le
\sum_{j\in GC} 2\int_{q_{z_j}}|\EE(\tilde u-w_{z_j})*\phi_\delta|^2\,dx
+2\int_{q_{z_j}} |\EE(w_{z_j}*\phi_\delta)-\phi_\delta*\E(u)|^2\,dx
\\\le
\sum_{j\in GC} \frac{c}{\delta^2}\int_{\tilde{q}_{z_j}}|\tilde u-w_{z_j}|^2\,dx
+
 c\left(\frac{\H^{n-1}(J_u\cap Q_{z_j})}{\delta^{n-1}}\right)^q
\hspace{-.5em}\int_{Q_{z_j}} |\E(u)|^2\,dx
\\\le
\sum_{j\in GC} \frac{c}{\delta^2}\int_{\tilde{q}_{z_j}}|\tilde u-w_{z_j}|^2\,dx
+ c\theta^q \int_{\Omega'}|\E(u)|^2\,dx
\end{multline}
thanks to~\eqref{estimdif}. 
{We now evaluate $\tilde u-w_{z_j}$ in $\tilde{q}_{z_j}$. If
$x\in \tilde{q}_{z_j}\setminus \cup_i \omega_{z_i}$, then $\tilde u(x)=w_{z_j}(x)=u(x)$ and the difference vanishes.
If $x\in \omega_{z_j}\setminus \cup_{i<j} \omega_{z_i}$ then $\tilde u(x)=w_{z_j}(x)=a_{z_j}(x)$ and again the difference vanishes. 
The remaining contributions are $a_{z_i}-a_{z_j}$ on the set $ \omega_{z_j}\cap \omega_{z_i}\setminus \cup_{k<i} \omega_{z_k}$, if $i<j$,
and $a_{z_i}-u$ on $\omega_{z_i}\setminus \cup_{k<i} \omega_{z_k}$ for $i>j$.} Hence
we can bound the integrals in the sum as follows:
\begin{multline*}
%\int_{q_{z_j}}|\EE((\tilde u-w_{z_j})*\phi_\delta)|^2\,dx\le
\frac{c}{\delta^2}\int_{\tilde q_{z_j}}|\tilde u-w_{z_j}|^2\,dx
= \frac{c}{\delta^2}\left(\sum_{i<j} |{\omega_{z_j}\cap\omega_{z_i}\setminus\cup_{k<i}\omega_{z_k}}|
{ \|a_{z_i}-a_{z_j}\|_{L^\infty(Q_{z_i}\cap Q_{z_j})}^2} \right.\\\left.
+ \sum_{i>j} \int_{\tilde q_{z_j}\cap(\omega_{z_i}\setminus \cup_{k<i}\omega_{z_k})}|u-a_{z_i}|^2 dx
\right).
\end{multline*}
The sums above involve at most $3^n-1$ terms corresponding
to the {\it good}  cubes that are neighbors with $q_{z_j}$. Thanks to Eq.~\eqref{eq:diffazz} and the fact that  the $Q_{z_j}$ are {\it good}, we have the
bound
\begin{multline*}
|{\omega_{z_j}\cap\omega_{z_i}\setminus\cup_{k<i}\omega_{z_k}}|\; \|a_{z_i}-a_{z_j}\|_{L^\infty(Q_{z_i}\cap Q_{z_j})}^2
 \\ \le
c\theta \delta^n
\frac{C}{\delta^{n-2}}\int_{Q_{z_i}\cup Q_{z_j}}\hspace{-1em}|\E(u)|^2dx
\le C\theta\delta^2\int_{Q_{z_i}\cup Q_{z_j}}\hspace{-1em}|\E(u)|^2dx.
\end{multline*}
Further, recalling Eq.~\eqref{eq:estimaz} and using, once again, the fact that $Q_{z_j}$ is a {\it good} cube,
\begin{multline*}
\int_{\tilde q_{z_j}\cap(\omega_{z_i}\setminus \cup_{k<i}\omega_{z_k})}|u-a_{z_i}|^2 dx
\le
|\tilde q_{z_j}\cap\omega_{z_i}|^{1/n} C\delta\int_{Q_{z_j}} |\E(u)|^2 dx
\\
 \blue{\le C(\delta\H^{n-1}(J_u\cap Q_{z_i}))^{1/n}\delta \int_{Q_{z_j}} |\E(u)|^2 dx}
 \le C\theta^{1/n}\delta^2 \int_{Q_{z_j}} |\E(u)|^2 dx.
\end{multline*}
{As a consequence,
\begin{equation}\label{L2estu-w}
\frac{c}{\delta^2}\int_{\tilde q_{z_j}}|\tilde u-w_{z_j}|^2\,dx
\le C\theta^{1/n}\int_{\tilde Q_{z_j}}|\E(u)|^2 dx
\end{equation}
where $\tilde Q_{z_j}$ is the cube $z_j+(-8\delta,8\delta)^n$, so
we deduce from~\eqref{eq:globalqzdiff} that
\begin{equation}\label{eq:errorstrain}
\int_{\Omega_g^\delta}|\EE(u_\e)-\phi_\delta*\E(u)|^2\,dx
\le c\theta^{q'} \int_{\Omega'}|\E(u)|^2\,dx
\end{equation}
where $q'=\min\{1/n,q\}$, for some constant $c>0$.
}
It follows easily from~\eqref{eq:errorstrain} that
\begin{equation}\label{eq:estimglobale}
\limsup_{\theta\to 0}\limsup_{\e \to 0} \int_{\Omega_g^\delta}(|\EEd(u_\e)|^2 + (\Dp u_\e)^2) dx\le \int_{\Omega'}(|\Ed(u)|^2+(\Trp\E(u))^2)\;  dx.
\end{equation}

\begin{remark}\label{l2conv}
Remark that we  have proved along the way that 
$$\limsup_{\theta\to 0}\limsup_{\delta\to 0} \|\tilde u -u\|_{L^2(\Omega')} =0$$
as could be easily checked using Eq.~\eqref{L2estu-w}, the control of $\omega_z$ by $\delta\H^{n-1}(J_u\cap Q_z)$ and the fact that there are a finite ($\delta$-independent) number of overlaps between the $Q_z$.
\end{remark}

\bigskip

We now address the union $\Omega_b^\delta$ of the {\it bad} cubes\footnote{{Which now includes
former ``good'' cubes $Q_z$ where $v_\e^\ell\equiv 0$}}.  There, we saw that, each {\it bad} cube has as  diameter  at most $2 \sqrt n \delta$, $v_\e^\ell\equiv 0$ and thus we only have to estimate from above
$$
\eta_\e\int_{\Omega_b^\delta}  \left\{(|\EEd(u_\e)|^2 + (\Dp u_\e)^2 )\right\}dx
\le 2\eta_\e |{\Omega_b^\delta}|c^2M^2/\delta^2.
$$
{Indeed, as  in Subsection \ref{secfirstlimsup},  $|\EE(u_\e)|\le c\|\tilde u\|_{L^\infty(\Omega')}/\delta
\le cM/\delta$.}

Choosing $\eta_\e=o(\delta)$ ensures that
\begin{equation}\label{eq:estimbad}
\limsup_{\e\to 0} \int_{\Omega_b^\delta}{(\eta_\e+(v_\e)^2)}  \left\{(|\EEd(u_\e)|^2 + (\Dp u_\e)^2) \right\}\ dx=0,
\end{equation}
{thanks to \eqref{eq:vol-v=0}.}

\bigskip

In view of Eqs.~\eqref{eq:estimvpart}, \eqref{eq:estimglobale}, \eqref{eq:estimbad}, 
we have proved that
\[
\limsup_{\theta\to 0}\limsup_{\e \to 0}\int_{\Omega'}{(\eta_\e+(v_\e)^2)}(|\EEd(u_\e)|^2 + (\Dp u_\e)^2)\; dx\le\int_{\Omega'}(|\Ed(u)|^2+(\Trp\E(u))^2)\;  dx.
\]

 The desired conclusion would be achieved  if we could show that 
\begin{equation}\label{eq:vanishdivm}
\limsup_{\e\to 0} \int_{\Omega'} (\Dm u_\e)^2 dx \le
 { \int_{\Omega'} (\Trm\E(u))^2 dx = }
\int_{\Omega'} |\Dm u|^2 dx
\end{equation}
({since} the negative part of divergence \blue{of $u$} has no singular part).

%As before in Eq.~\eqref{eq:estimglobale}, we would get, choosing this time $T=\text{tr}^-\; (p=1)$,
{As before, thanks to~\eqref{eq:errorstrain} we have}
\begin{equation*}
\limsup_{\theta\to 0}\limsup_{\e \to 0} \int_{\Omega_g^\delta} (\Dm u_\e)^2dx \le \limsup_{\e \to 0}\int_{\Omega'}|\Dm u|^2dx.
\end{equation*}
We are thus only concerned with $\int_{\Omega_b^\delta} (\Dm u_\e)^2dx$. That quantity could immediately be seen to go to $0$ if  we had that $\tilde u\equiv u$ on $\Omega_b^\delta$ because $\lim_{\e \to 0} |\Omega_b^\delta|=0$ while $\Dm u\in L^2(\Omega)$. Unfortunately such is not the case because the small  sets $\omega_{z_i}$ are only included in $\tilde q_{z_i}$.

A different route might consist in decomposing
\begin{multline*}
\Div u_\e = \phi_\delta*\Div u + \Div [\phi_\delta*(\tilde u-u)]=\phi_\delta*(\Dp u -\Dm u)+ \Div [\phi_\delta*(\tilde u-u)]
\ge\\ -\phi_\delta*\Dm u -| \Div [\phi_\delta*(\tilde u-u)] |,
\end{multline*}
so that
\[
\Dm u_\e \le (\phi_\delta*\Dm u) + | \Div [\phi_\delta*(\tilde u-u)] |.
\]
Since $\Dm u\in L^2(\Omega')$, it is immediate that
{\begin{equation}\label{div:nomodif}
\limsup_{\e\to 0} \int_{\Omega_b^\delta} (\phi_\delta*\Dm u)^2 dx = 0. %\int_{\Omega'} (\Dm u)^2 dx.
\end{equation}}
%
%THIS SPLITTING IS NOT {\it good}   SINCE it is not clear that $\tilde u$
%is close enough to $u$...

It would remain to prove that
$\| \Div [\phi_\delta*(\tilde u-u)] \|_{L^2({\Omega_b^\delta})}$  yields a negligible contribution in the limit. Recall that
\[
\tilde u (x)-u(x) = \begin{cases}
0 & \textup{ if } x\in   %%\bigcup_{j\in J} \tilde q_{z_j} 
\Omega' \setminus \bigcup_{j\in \GC} \omega_{z_j} \\
a_{z_j}-u(x) & \textup{ if } x \in \omega_{z_j} \setminus \bigcup_{i<j} \omega_{z_i}.
\end{cases}
\]
so that
\[
|\tilde u(x)-u(x)|\le 2M\sum_{j\in \GC}\chi_{\omega_{z_j}}.
\]
Hence, observing that $\phi_\delta*\chi_{\omega_{z_j}}$ has 
support in $Q_{z_j}$ and that the $Q_{z}$ have finite overlap,
\begin{align*}
\int_{{\Omega_b^\delta}} |\Div [\phi_\delta*(\tilde u-u)] |^2
& \le 4M^2 \int_{{\Omega_b^\delta}} [{|\nabla\phi_\delta|}*\sum_{j\in \GC}\chi_{\omega_{z_j}}]^2dx
\\
& \le CM^2 \sum_{j\in \GC}\int_{{\Omega_b^\delta}} [{|\nabla\phi_\delta|}*\chi_{\omega_{z_j}}]^2dx
\\
& \le C\frac{M^2}{\delta^2}\sum_{j\in \GC} |\omega_{z_j}\cap (\Omega_b^\delta+B(0,\tfrac{\delta}{2}))|\\
&\le C\frac{M^2}{\delta^2}  \delta\sum_{j\in \GC} \H^{n-1}(J_u\cap Q_{z_j})
\\ & \le C\frac{M^2}{\delta} \H^{n-1}(J_u{\setminus\Gamma})\end{align*}
so that 
\[
\| \Div [\phi_\delta*(\tilde u-u)] \|_{L^2({\Omega_b^\delta})}
\le C\frac{M}{\sqrt{\delta}}{\theta}.
\]
That estimate also fails.
We do not know at present how to circumvent this difficulty for an arbitrary dimension $n$.

\subsection{The two-dimensional case}\label{sec:2-d-case}
In two dimension, we make  use  of Theorem~2.1~in~\cite{CFI}, which we specialize here
to  $p=2$ and restate in a simpler form.
\begin{theorem}\label{thm:CFI}
There exist $\zeta>0$ and $\check c>0$ such that for $u\in SBD^2(B_{2r})$
with $\H^1(J_u)\le 2\zeta r$, there exists $R\in (r,2r)$ and
$\check u\in SBD^2(B_{2r})\cap \blue{H^1}(B_R;\R^2)$ such that\begin{itemize}
\item[(i)] $\blue{\H^1}(J_{u}\cap \partial B_R)=0$;
\vskip.2cm

\item[(ii)] $
\int_{B_R}|\E(\check u)|^q dx \le \check c \int_{B_R}|\E(u)|^q dx$ for
$q\in [1,2]$;
\vskip.2cm

\item[(iii)] $\|u-\check u\|_{L^1(B_R;\R^2)} \le\check c R|\ED u|(B_R)$;
\vskip.2cm

\item[(iv)] $u=\check u$ on $B_{2r}\setminus B_R$, 
$\blue{\H^1}(J_{\check u}\cap \partial B_R)=0$;
\vskip.2cm

\item[(v)] if $u\in L^\infty(B_{2r},\R^2)$ then
{$\|\check u\|_{L^\infty} \le \|u\|_{L^\infty}$.}
\end{itemize}
\end{theorem}

Thanks to this result, we can slightly amend the construction
in Section~\ref{sec:GenProof} and conclude. 
%The issue was to connect the ``good'' part where we could modify slightly $u$ to obtain a global control on $\int |\E(u)|^2dx$ with the rest of the domain, maintaining a global $L^2$ bound on the negative part of the divergence. 
The theorem allows to build a safety zone
where $u$ has no jump at all between $\Omega_g^\delta$ and $\Omega_b^\delta$.

The proof follows the same lines until Eq.~\eqref{eq:estimvpart}.
Then, one considers the ``boundary'' {\it good}   cubes, which are {\it good}   cubes
$q_z$ on which $v_\e\equiv 0$. By construction, $v_\e$ is zero in 
a $16\sqrt{2}\delta$-neighborhood of
$\Sigma:=\Gamma\cup \{ z: Q_z \textup{ is bad}, Q_z\cap \Gamma=\emptyset\}$,
and positive elsewhere.
Hence, with the additional help of Eq.~\eqref{badcubes}, we conclude that there are at most $C(\H^1(J_u)+\theta)/\delta\le C'/\delta$ ($C'$ is a different constant) such cubes, so
that they cover a surface area of order $\delta$.
We call $\BGC$ a numeration of these {\it good}   cubes $\{Q_{\tilde z_i}\}_{i\in \BGC}$.
Assuming
$\theta < 2\zeta$ where $\zeta$ is the constant in Theorem~\ref{thm:CFI},
we build recursively a function $u_i$ as follows: we let $u_0=u$
and for each $i\in\BGC$,
in $B(\tilde z_i,4\sqrt{2}\delta)$,
we find $\check u_i$ such that $\check u_i=u_{i-1}$
near $\partial B(\check z_i, 4\sqrt{2}\delta)$,
$\check u_i\in H^1(\check q_{\check z_i})$,
\[
\int_{B(\check z_i,4\sqrt{2}\delta)}|\E(\check u_i)|^2dx \le
\check c\int_{B(\check z_i,4\sqrt{2}\delta)}|\E(u_{i-1})|^2dx ,
\]
and
\[
\|u_{i-1}-\check u_i\|_{L^1(B(\check{z}_i,4\sqrt{2}\delta);\R^2)} \le\check c \delta|\ED u_{i-1}|(B(\check z_i,4\sqrt{2}\delta)).
\]
We then let $u_i=\check u_i$ in $B(\check z_i,4\sqrt{2}\delta)$ and $u_i=u_{i-1}$
in the rest of the domain $\Omega'$. We call $I:=\#(\BGC)$. Note that the balls $B(\check z_j,4\sqrt{2}\delta)$
overlap a finite number of times, hence, possibly changing $\check c$, one
still has for $i\in \BGC$
\[
\int_{B(\check z_i,4\sqrt{2}\delta)}|\E(u_I)|^2dx \le
\check c\int_{B(\check z_i,4\sqrt{2}\delta)}|\E(u)|^2dx ,
\]
and
\[
\|u- u_I\|_{L^1(B(\check{z}_i,4\sqrt{2}\delta);\R^2)} \le\check c \delta|\ED u|(B(\tilde z_i,4\sqrt{2}\delta)).
\]
It follows that 
\begin{equation}\label{estl1uiu}
\|u-u_I\|_{L^1(\Omega')} \le \check c \delta
\end{equation}
and, setting $J^\delta := \bigcup_{i\in\BGC} B(\check z_i,4\sqrt{2}\delta)$,
\[
\int_{J^\delta}|\E(u_I)|^2dx \le\check c \int_{J^\delta}|\E(u)|^2dx.
\]
As $|J^\delta|\le C \delta$, this quantity can be made arbitrarily
small: we deduce that
\begin{equation}\label{estimui}
\int_{\Omega'}|\E(u_I)|^2dx \le \int_{\Omega'} |\E(u)|^2 + {\theta} %C\delta.
\end{equation}
{for  $\delta$ small enough.}
We remark that, if  $\check J^\delta:=
\bigcup_{i\in\BGC}  q_{\check z_i}$, then by construction $u_I\in H^1(\check J^\delta)$ and also that
\begin{equation}\label{difuiu}
\{u_I\ne u\}\subset J^\delta.
\end{equation}

Now we start the very construction of Subsection~\ref{sec:GenProof}
after Eq.~\eqref{eq:estimvpart}  with $u$ replaced with $u_I$. 
The only difference is that in the {\it good}   cubes $q_{\tilde z_i}$,
\blue{ $i\in \BGC$,}  (for which $u\in H^1(q_{\tilde z_i})$)
we are at liberty to set $\omega_{\tilde z_i}=\emptyset$. Thanks to Eq.~\eqref{estimui} and to item (v) in Theorem \ref{thm:CFI}, the estimates in Eqs.~\eqref{eq:estimglobale}, \eqref{eq:estimbad} still hold true and, further, thanks to Eq.~\eqref{difuiu},
\begin{equation}\label{diftildeuiu}
\{\tilde u_I\ne u\}\subset J_\delta.
\end{equation}

Finally, we have to estimate  $\int_\Omega (\Dm u_\e)^2dx$.
As before, { from Eq.~\eqref{eq:errorstrain}} we would get
\begin{equation}\label{div-og}
\limsup_{\theta\to 0}\limsup_{\e \to 0} \int_{\Omega_g^\delta} (\Dm u_\e)^2dx \le \limsup_{\e \to 0}\int_{\Omega'}|\Dm(u_I)|^2dx =\int_{\Omega'}|\Dm(u)|^2dx
\end{equation}
where the last equality results from  Eqs.~\eqref{estimui},\eqref{diftildeuiu}.

If $x\in\Omega\setminus\Omega_g^\delta$, $x$ is either in a {\it bad} cube $q_{z_j}$
and $\dist(x,\Sigma)\le \sqrt{2}\delta$, or $x$ is in a cube $Q_z$
which intersects $\Gamma$ so that $\dist(x,\Sigma)\le 8\sqrt{2}\delta$.
Hence $\Omega\setminus\Omega_g^\delta\subset\{\dist(\cdot,\Sigma)\le 8\sqrt{2}\delta\}$.

Consider $q_{z_i}$ a {\it good}   cube, with
$\omega_i\neq\emptyset$. This means that $q_{z_i}\not\subset
\{\dist(\cdot,\Sigma)\le 16\sqrt{2}\delta\}$ so that
$\dist(q_{z_i},\Sigma)\ge 14\sqrt{2}\delta$.
It follows that $\dist(q_{z_i},\Omega\setminus\Omega_g^\delta)\ge 6\sqrt{2}\delta$,
and since $\omega_i\subset\tilde q_{z_i}$,
$\dist(\omega_i,\Omega\setminus\Omega_g^\delta)\ge 5\sqrt{2}\delta$.

We conclude that $\dist(\{u_I\neq \tilde u_I\},\Omega\setminus \Omega_g^\delta)\ge 5\sqrt{2}\delta$.
It means that there is a strip of width $5\sqrt{2}\delta$ where
$\tilde u_I = u_I$ is $H^1$ along the boundary between $\Omega\setminus \Omega_g^\delta$
and $\Omega_g^\delta$. Hence everywhere in $\Omega\setminus \Omega_g^\delta$, $u=u_I$ so that
$\Dm u_\e (x)\le (\phi_\delta*\Dm u)$ and it follows
that
\begin{equation}\label{div-ob}
\limsup_{\e\to 0}\int_{\Omega\setminus\Omega_g^\delta} (\Dm u_\e)^2dx \le
 \limsup_{\e\to 0}\int_{\{\dist(\cdot,\Sigma)\le 9\sqrt{2}\delta\}} (\phi_\delta*\Dm u)^2dx = 0.
\end{equation}
Recalling Eqs.~\eqref{div-og},\eqref{div-ob}, we conclude that, this time Eq.~\eqref{eq:vanishdivm} holds.

\vskip1cm

Finally note that
 $u_\e\to u$ 
strongly in $L^2(\Omega';\R^2)$. Indeed, recalling Eq.~\eqref{estl1uiu} and item (v) in Theorem \ref{thm:CFI}, 
$\|u_I-u\|_{L^2(\Omega')}\le C\delta$ while, by Remark \ref{l2conv}, $\|\tilde u_I-u_I\|_{L^2(\Omega')}\to 0$.

\section*{Acknowledgements}
\noindent 
G.A.F.'s research has been  supported {in part} by the National Science Fundation Grant DMS-1615839, {the research of S.C. was partially supported 
by the Deutsche Forschungsgemeinschaft through the Sonderforschungsbereich 1060
{\sl ``The mathematics of emergent effects''}.}

The authors warmly acknowledge the support and hospitality
of the Matematisches Forschungsinstitut Oberwolfach  where this
research was initiated (workshop \#0709) and completed (worskhop \#1729).
They also wish to thank Giovanni Lancioni for  providing the images of Fig.~\ref{fig:LRC}
as well as appropriate references, Blaise Bourdin for  further references
and also for very enlightening discussions, and Michael Ortiz for 
very helpful remarks on the models.

\appendix
\section{Affine approximation of bounded functions}

 \begin{lemma}\label{lemmalinfty}
  Let $R>r>0$, $Q_R:=x_*+(-R/2,R/2)^n$, $Q_r:=x_*+(-r/2,r/2)^n$, $\omega\subset Q_R$ such that
  $|\omega|\le (R-r)^n/2^{n+1}$, $A:\R^n\to\R^n$ affine, $u\in L^\infty(Q_R;\R^n)$. Then there is $a:\R^n\to\R^n$ affine 
  such that
  \begin{equation}
   \|a\|_{L^\infty(Q_r;\R^n)}\le \|u\|_{L^\infty(Q_R;\R^n)} 
  \end{equation}
  and
  \begin{equation}
   \|u-a\|_{L^p(Q_R\setminus\omega;\R^n)}\le c \|u-A\|_{L^p(Q_R\setminus\omega;\R^n)}.
  \end{equation}
  The constant $c\ge 1$ depends on $p$, $n$ and $r/R$. If $\EE(A)=0$, then $\EE(a)=0$.
 \end{lemma}
\begin{figure}[htb]
 \begin{center}
\includegraphics[width=6cm]{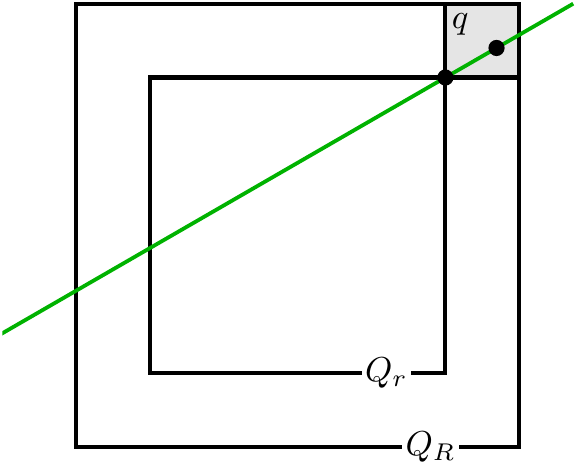}
 \end{center}
 \caption{Sketch of the geometry in the proof of Lemma \ref{lemmalinfty}. The line is $t\mapsto v+t(y-v)$, {the two dots are the} points $v$ and $y$.}
\label{figlinfty}
 \end{figure}
\begin{proof}
We can assume $x_*=0$ and    $\|A\|_{L^\infty(Q_r;\R^n)}> \|u\|_{L^\infty(Q_R;\R^n)}$, otherwise $a=A$ will do. 

 The function $x\mapsto |A(x)|$ is convex, therefore there is a vertex $v$ of $Q_r$ such that
 $|A(v)|= \|A\|_{L^\infty(Q_r;\R^n)}$. 
 Consider the cube
 \begin{equation*}
  q:=\frac{R+r}{2r} v + (-\frac{R-r}4,\frac{R-r}4)^n.
 \end{equation*}
We check that $q\subset Q_R\setminus Q_r$ and that $v$ is a vertex of $q$ as well (see Figure \ref{figlinfty}).
Further, for any $y\in q$ there is $t_*<0$ such that
$v+t(y-v)\in Q_r$ for all $t\in[t_*,0)$ (it suffices to check this componentwise, separating the cases $v_i=r/2$ and $v_i=-r/2$).
The function
\begin{equation*}
 t\mapsto f_y(t):= |A ( v + t (y-v))|
\end{equation*}
is convex, and obeys 
\begin{equation*}
 f_y(t_*)\le \|A\|_{L^\infty(Q_r)} = |A(v)|=f_y(0)
\end{equation*}
therefore $|A(y)|=f_y(1)\ge  |A(v)|$, which implies $|A(y)-u(y)|\ge |A(v)|-\|u\|_{L^\infty(Q_R)}>0$ for any $y\in q$.
 
Since $|\omega|\le \frac12 |q|$, we obtain
\begin{equation*}
 \frac12 |q| (|A(v)|-\|u\|_{L^\infty(Q_R)})^p\le \int_{Q_R\setminus\omega} |A-u|^pdx.
\end{equation*}
We define
\begin{equation*}
 a:=\frac{\|u\|_{L^\infty(Q_R)}}{|A(v)|}A
\end{equation*}
and estimate $\|a\|_{L^\infty(Q_r)}=\|u\|_{L^\infty(Q_R)}$ and
\begin{align*}
 \int_{Q_r} |A-a|^pdx &=\int_{Q_r} \left| \frac{|A(v)|-\|u\|_{L^\infty(Q_R)}}{|A(v)|} A\right|^pdx\\
& \le \int_{Q_r} (|A(v)|-\|u\|_{L^\infty(Q_R)})^pdx \le \frac{2|Q_r|}{|q|} \int_{Q_R\setminus\omega} |A-u|^pdx.
\end{align*}
Finally, since $A-a$ is affine we have, for some $c\ge 1$\footnote{It is
actually possible to prove that one can take $c=1$.}:
\begin{equation*}
 \int_{Q_R} |A-a|^pdx \le c\frac{R^{n+p}}{r^{n+p}} \int_{Q_r} |A-a|^pdx \le c\frac{2^{n+1}R^{n+p}}{r^p (R-r)^n} \int_{Q_R\setminus\omega} |A-u|^pdx.
\end{equation*}
A triangular inequality concludes the proof.
\end{proof}

\end{document}